\documentclass[12pt]{amsart}

\usepackage{amssymb}
\usepackage{amsmath}

\evensidemargin\oddsidemargin

\begin{document}

\baselineskip=18pt
\setcounter{page}{1}

\newtheorem{LEM}{Lemma\!\!}
\newtheorem{THA}{Theorem A\!\!}
\newtheorem{THB}{Theorem B\!\!}
\newtheorem{CORO}{Corollary\!\!}
\newtheorem{CONJ}{Conjecture \!\!}
\newtheorem{CONC}{Conjecture C\!\!}
\newtheorem{COND}{Conjecture D\!\!}
\newtheorem{REM}{Remark\!\!}

\renewcommand{\theLEM}{}
\renewcommand{\theTHA}{}
\renewcommand{\theTHB}{}
\renewcommand{\theCORO}{}
\renewcommand{\theCONC}{}
\renewcommand{\theCOND}{}
\renewcommand{\theCONJ}{}
\renewcommand{\theREM}{}

\newcommand{\eqnsection}{
\renewcommand{\theequation}{\thesection.\arabic{equation}}
    \makeatletter
    \csname  @addtoreset\endcsname{equation}{section}
    \makeatother}
\eqnsection

\def\AA{\mathcal{A}}
\def\BB{\mathcal{B}}
\def\CC{\mathcal{C}}
\def\DD{\mathcal{D}}
\def\d{{\mathrm d}}
\def\ddd{{\mathbb D}}
\def\e{{\mathbb E}}
\def\eps{\varepsilon}
\def\FF{\mathcal{F}}
\def\GG{\mathcal{G}}
\def\LL{\mathcal{L}}
\def\LLA{\LL}
\def\hA{\hat A}
\def\hX{\hat X}
\def\hZ{\hat Z}
\def\k{{\mathcal K}}
\def\p{{\mathbb P}}
\def\r{{\mathbb R}}

\def\lacc{\left\{}
\def\lcr{\left[}
\def\lpa{\left(}
\def\lva{\left|}
\def\racc{\right\}}
\def\rcr{\right]}
\def\rpa{\right)}
\def\rva{\right|}

\def\Un{{\bf 1}}
\def\ee{\mathrm{e}}
\def\qed{\hfill$\square$}
\def\elaw{\stackrel{d}{=}}

\title[Burgers equation with stable initial data]
      {On the Hausdorff dimension of regular points of inviscid Burgers equation with stable initial data}

\author[Thomas Simon]{Thomas Simon}

\address{Equipe d'Analyse et Probabilit\'es, Universit\'e d'Evry-Val d'Essonne, Boulevard Fran\c{c}ois Mitterrand, F-91025 Evry Cedex. {\em E-mail address}: {\tt tsimon@univ-evry.fr}}



\begin{abstract} Consider an inviscid Burgers equation whose initial data is a L\'evy $\alpha-$stable process $Z$ with $\alpha > 1$. We show that when $Z$ has positive jumps, the Hausdorff dimension of the set of Lagrangian regular points associated with the equation is strictly smaller than $1/\alpha,$ as soon as $\alpha$ is close to 1. This gives a negative answer to a conjecture of Janicki and Woyczynski \cite{JW}. Along the way, we contradict a recent conjecture of Z. Shi about the lower tails of integrated stable processes.
\end{abstract}

\maketitle

\section{Introduction and statement of the results}
Since the seminal paper \cite{Si1}, statistical properties of the Burgers equation 
\begin{equation}
\label{Burg}
\partial_t u \; +\; u\partial_x u \; =\; \nu \partial_{xx} u, \qquad \nu >0
\end{equation}
with initial condition $u_0(x) := u(0,x) = X_x$ where $\lacc X_x, \, x\in\r\racc$ is a given random process, have given rise to intensive research. Eventhough (\ref{Burg}) is a much simplified version of the Navier-Stokes equation, it is still relevant in physics as a model equation for e.g. shock waves in hydrodynamics, in order to describe turbulence phenomena not covered by the classical Korteweg-de Vries equation. From the mathematical point of view, a nice feature of (\ref{Burg}) is the possibility to solve it explicitly through the change of variable $u = -\partial_x \psi$, which is known as the Hopf-Cole substitution: one has 
\begin{equation}
\label{HC}
\psi(t,x)\; =\; 2\nu \log\lcr (4\pi\nu t)^{-1/2} \int_\r \exp \lcr \frac{1}{2\nu} \lpa \psi_0(a) - \frac{(x-a)^2}{2t}\rpa\rcr da \rcr
\end{equation}
where $\psi_0$ is the initial potential given by $u_0 = -\partial_x \psi_0$. 

The {\em inviscid} Burgers equation is a simplified version of (\ref{Burg}) where the viscosity parameter $\nu = 0,$ and its so-called Hopf-Cole solution is obtained from (\ref{HC}) in letting $\nu\to 0$. By Laplace approximation, it takes a particularly nice form:
\begin{equation}
\label{HCbis}
\psi(t,x)\; =\; \sup_{a\in\r} \lacc \psi_0(a) - \frac{(x-a)^2}{2t}\racc,
\end{equation}
which is well-defined provided the initial potential satisfies $\psi_0(a) = {\rm o} (a^2)$ when $\vert a\vert\to  +\infty.$ We refer e.g. to the monograph \cite{Wo} for the above facts, and much more, concerning Burgers equation. 

In this paper, we are interested in the inviscid Burgers equation whose initial data is a two-sided L\'evy $\alpha-$stable process. More precisely, we suppose that the initial condition $X$ is defined followingly:
\begin{equation}\label{XX}
X_x \; =\; \lacc\begin{array}{cl} Z_x & \mbox{if $x\ge 0$}\\
                                                   -Z'_{-x} & \mbox{if $x\le 0$}\end{array}\right.\end{equation}
where $Z =\lacc Z_x, \, x\ge 0\racc$ is a L\'evy $\alpha-$stable process and $Z'$ an independent copy of $Z$. Referring to Chapter VIII in \cite{Ber} for more details, let us recall that $Z$ is a real process with stationary and independent increments, which is $(1/\alpha)-$self-similar:
\begin{equation}
\label{SS}
\{ Z_{kx}, \; x\ge 0\}\; \elaw\;\{ k^{1/\alpha}X_x, \; x\ge 0\}\end{equation}
for all $k >0$. This property forces the stability index $\alpha$ to be in $(0,2],$ and $Z$ is Brownian motion (up to a scaling parameter) in the case $\alpha =2.$ Its L\'evy-Khintchine exponent $\Psi(\lambda) = - \log \e [\ee^{i\lambda Z_1}]$ is given by
$$\Psi(\lambda) \; =\;\kappa\vert\lambda\vert^\alpha(1-{\rm i}\beta{\rm sgn}(\lambda)\tan(\pi\alpha/2)), \quad \lambda\in\r,$$
where $\kappa > 0$ is the scaling parameter and $\beta\in [-1,1]$ is the skewness parameter. Without loss of generality, in the following we will suppose $\kappa \equiv 1.$ The positivity parameter $\rho\; =\; \p\lcr Z_1 > 0\rcr$ takes its values in $[1-1/\alpha, 1/\alpha]$ when $\alpha > 1$ and in $[0,1]$ when $\alpha <1$. When $\alpha >1,$ the value $\rho = 1- 1/\alpha$ corresponds to the spectrally positive situation ($\beta = 1$ and $Z$ has no negative jumps) and the value 
$\rho = 1/\alpha$ to the spectrally negative situation ($\beta =-1$ and $Z$ has no positive jumps). When $\alpha < 1$ and $\rho = 0$ (resp. $\alpha < 1$ and $\rho =1$), $Z$ is the negative of a subordinator and has no positive jumps (resp. a subordinator and has no negative jumps). When $\alpha =1$, then $\rho\in (0,1)$ and $Z$ has jumps in both directions. 

The law of the iterated logarithm for $Z$ at infinity - see Theorem VIII.5 in \cite{Ber} - entails 
$$\limsup_{x\to +\infty} x^{-\kappa}Z_x\, =\, 0\;\; \mbox{or}\;\; +\infty\quad\mbox{according as}\quad \kappa > 1/\alpha\;\; \mbox{or}\;\;\kappa \le 1/\alpha.$$
Hence, since our initial potential is given up to some meaningless additive constant by 
$$\psi_0(x)\; = \;\int_0^x X_t\, \d t, \qquad x\in\r,$$  
the growth condition $\psi_0(a) = {\rm o} (a^2)$ at infinity assigns the restriction 
$$\alpha\in (1,2]$$ 
on the stability parameter, which will be supposed henceforth unless explicitly stated. 

Derivating (\ref{HCbis}) with respect to $x$ yields readily the following formula for the Hopf-Cole solution of (\ref{Burg}):
$$u(t,x) \; = \; \frac{x - a(t,x)}{t}$$
where $a(t,x)$ is the largest point attaining the maximum in (\ref{HCbis}), in other words:
$$a(t,x)\; =\; \max\lacc s\in\r, \; C'_t(s)\le x t^{-1}\racc$$
where $C'_t$ stands for the right-derivative of $C_t,$ which is the convex hull of the function
$$x\;\mapsto\; \int_0^x (X_u + u t^{-1}) \d u.$$
The so-called {\em Lagrangian regular points} of (\ref{Burg}) are the points where the above function coincides with its convex hull. Notice that this time-dependent set 
$\LLA_t$ can be described in terms of the function $a(t,x)$:
$$\LLA_t\; =\; \lacc a(x,t), \, x\in\r\;\;\mbox{and}\;\; a(x-,t) = a(x,t)\racc.$$ 
From the physical point of view, $\LLA_t$ is the set of particles which have not been shocked up to time $t$ by the turbulence governed by Equation (\ref{Burg}) - see \cite{Wo} - and for this reason there has been some interest over the years in describing the structure of the sets $\LLA_t$, especially in studying their fractal properties. In \cite{JW}, the authors raised the following 

\begin{CONJ}[Janicki and Woyczynski] For every $t >0,$ one has
$${\rm Dim}_H \, \LLA_t\; = \; 1/\alpha\quad\mbox{a.s.}$$
\end{CONJ} 

It is easy to see that the Hausdorff dimension of $\LLA_t$ does not, indeed, depend on $t$ in this model. Namely, from the self-similarity of $Z$ one can show - see \cite{JW} p. 285 - that
$$u(t,x)\;\elaw\;t^{1/(\alpha-1)}u(1, xt^{\alpha/(\alpha-1)}),$$
which entails that a.s. ${\rm Dim}_H \, \LLA_t = {\rm Dim}_H \, \LLA_1$ for every $t >0$. In the following, we will be therefore interested in the set $\LLA_1$ only, which we denote by $\LLA$ for the sake of simplicity. 

Notice that the above conjecture had been previously solved (without complete proof) by Sinai \cite{Si1} in the Brownian case $\alpha =2$, and that Bertoin \cite{Ber1} showed then rigourously the result in the general case $\alpha\in (1,2]$ with no positive jumps, as a consequence of the remarkable fact that the process $x\mapsto a(x,1)$ is a subordinator close to the$(1/\alpha)$-stable one - see Theorem 2 in \cite{Ber1}. Nevertheless, this latter property is no more true in the non spectrally negative framework - see the conclusion of \cite{Ber1}, and the structural study of (\ref{Burg}) seems to require an entirely different methodology when there are positive jumps. To this end, let us also mention an attempt made in \cite{CD2} with the concept of statistical solution (which is different from the Hopf-Cole solution).  

In this paper we show that Janicki and Woyczynski's conjecture is false in general when there are positive jumps:

\begin{THA} Suppose that $Z$ has positive jumps. Then there exists $\alpha_0 > 1$ such that
$${\rm Dim}_H \, \LLA\; < \; 1/\alpha\quad\mbox{a.s.}$$
for every $\alpha\in (1,\alpha_0).$
\end{THA}

Our main argument comes from a recent paper by Molchan and Kholkhov \cite{MK}, who were interested in the sets $\LL_t$ of Lagrangian regular points of (\ref{Burg}) when the initial data is a two-sided fractional Brownian motion $W$. In \cite{MK}, they proved that an a.s. upper bound on ${\rm Dim}_H \,\LL_t$ - which is also independent of $t > 0$ by the self-similarity of $W$ - follows from a lower bound on the exponent $\kappa$ appearing in the estimate
$$\p\lcr \int_0^t {\hat W}_s \d s < \eps \; + \; t^2,\;\forall\, t\in [-1,1] \rcr\;\le\;\eps^{\kappa}, \qquad \eps \to 0,$$
where ${\hat W} = -W \elaw W$ is the dual process of $W$. We remark
that this argument transfers to the L\'evy stable case without much
difficulty, and is actually even simpler because of the independence
and stationarity of the increments of $X$. This will be done in
Section 3. Before this, we will have to prove a crucial estimate on
the first-passage time of the integrated stable process. Let us fix
first some more notation, introducing 
$$A_t \; =\; \int_0^t Z_s \,\d s, \qquad t\ge 0,$$
the integral of a L\'evy $\alpha-$stable process $Z$ - with the same notations as above, but this time for every $\alpha\in (0,2]$ - and $T = \inf\{ t >0, \, A_t = 1\}$ the first-passage time of $A$ across 1.

\begin{THB} Suppose that $Z$ has negative jumps and is not the negative of a subordinator. With the above notations, for every $\alpha_0 > 0$ there exists a constant $\kappa > 0$ depending only on $\rho$ such that
$$\liminf_{t\to\infty} t^\kappa \p[ T > t]\; =\; 0$$
for every $\alpha \in [\alpha_0, 2].$  
\end{THB}

This result, which will be proved in Section 2, is interesting in its own right because it contradicts another conjecture whose solution had been announced (with a hidden error) by the author in \cite{ThS1}, and whose statement was the following:

\begin{CONJ}[Shi] Suppose that $Z$ is symmetric and that $\alpha >1$. Then with the above notations,
\begin{equation}
\label{sh}
\p[ T > t]\; =\; t^{-(\alpha -1)/2\alpha +o(1)}, \quad t\to\infty.
\end{equation}
\end{CONJ}
Indeed, since in Theorem B the exponent $\kappa$ is positive independent of $\alpha$, we see that the above conjecture is contradicted when $\alpha$ is close enough to 1, as in Theorem A. Notice that in a recent paper \cite{ThS} the author proved (\ref{sh}) in the case where $Z$ has no negative jumps and $\alpha >1$. The idea - which had been originally given by Z. Shi - consisted in time-changing the process $A$ through $\tau$ the inverse local time of $Z$ and considering the fluctuations of the L\'evy stable process $A_\tau$. We have been thinking for a long time that this method would be also successful when $Z$ has negative jumps, because $A_\tau$ is a L\'evy symmetric $(\alpha-1)/(\alpha +1)$-stable process {\em whatever} the value of $\rho$ should be - see Lemma 1 in \cite{ThS}, which allows in particular to obtain a general upper bound 
$$\p[ T > t]\; \le\; \k\, t^{-(\alpha -1)/2\alpha}, \quad t\to\infty
$$
for every value of $\rho\in [1-1/\alpha, 1/\alpha]$ and some finite constant $\k$ - see Theorem A in \cite{ThS}. It now appears that this subordination method is too crude when there are negative jumps, and yields only an upper bound which is not optimal, at least when $\alpha$ is close to 1. In this paper, we will obtain a uniform upper bound by discretization through exponential time-change combined with FKG-type inequalities, all of which we learned from Caravenna and Deuschel in the genesis of their paper \cite{CD1}. Let us stress that these arguments are also by far non optimal. Nevertheless they are quite robust and, since they involve eventually only fixed upper tails of $Z$ and $A$ which can be bounded independently of $\alpha$, this method makes it at least possible to contradict both Janicki-Woyczynski and Shi's conjectures when $\alpha$ is close to 1. Actually, we believe that these conjectures are false for all values of $\alpha$ except when there are no positive jumps (for the first) or no negative jumps (for the second), and in the fourth and last section of the paper we will state two other conjectures for the values of ${\rm Dim}_H \LL$ and the critical exponent in (\ref{sh}), in a general non completely asymmetric framework. 

\section{Proof of Theorem B}

We begin with the case $\alpha > 1,$ and we will actually obtain a
slightly stronger result which is crucial for Theorem A. Fix
$\rho = \p [Z_1 > 0] > 0$ once and for all, and set $\gamma = (\alpha -1)/\alpha > 0$. We will show that there exists $\kappa > 0$ independent of $\alpha$ such that 
\begin{equation}
\label{main}
\liminf_{t\to +\infty} t^\kappa \p\lcr A_s < 1 + s + t^{-\gamma} s^2, \; \;\forall\; s\le t\rcr\; =\; 0,
\end{equation}
which readily entails Theorem B by comparison. We first define an exponential subsequence of times, in considering the events
$$\AA_n\; =\; \lacc A_{2^m} < 1 + 2^m + 4^{m-n\gamma}, \; m = 0 \ldots 2n\racc, \quad n \ge 0.$$
Clearly, it is enough to prove that there exists $\kappa > 0$ independent of $\alpha$ such that $\p [ \AA_n ]\; \le\; \ee^{-\kappa n}$ for $n$ sufficiently large. We will obtain slightly more, in showing that 
\begin{equation}
\label{main2}
\p [ \BB_n ]\; \le\; \ee^{-\kappa n}
\end{equation}
for some $\kappa > 0$ independent of $\alpha$ and $n$ large enough, where 
$$\BB_n\; =\; \lacc A_{2^m} < 1 + 2^m + 4^{m-n\gamma}, \; m = 0 \ldots n\racc, \quad n \ge 0.$$
To do so, consider the events $\CC_k \; =\; \{ Z_{2^k} > - 2^{k/\alpha}, \; A_{2^k} > -2^{k(1+1/\alpha)}\}$ for every $k \ge 0,$ and the family of random times defined recursively by
$$\sigma_0\, =\, 0\quad\mbox{and}\quad \sigma_n\, = \, \inf\{ k > \sigma_{n-1}\; / \; \CC_k \; \mbox{occurs}\}, \quad n\ge 1.$$    
If $\{\FF_t, \, t\ge 0\}$ stands for the completed $\sigma-$field generated by $\{ Z_s, \, s\le t\},$ then $2^{\sigma_n}$ is a $\FF_t$-stopping time for every $n\ge 0.$ Denoting
henceforth by $[x]$ the integer part of any real number $x,$ we now state a crucial lemma which is mainly borrowed from \cite{CD1}: 

\begin{LEM} For every $\alpha_0 >0,$ there exists $\delta, c > 0$ independent of $\alpha\in [\alpha_0, 2]$ such that
$$\p\lcr \sigma_{[\delta n]} \, < \, n\; \vert \; \BB_n\rcr\; \ge \; c$$
for all $n$ sufficiently large.
\end{LEM}

Taking this lemma for granted and setting $K = c^{-1} < +\infty$ we see that
\begin{eqnarray*}
\p\lcr \BB_n\rcr & \le & K \p\lcr \sigma_{[\delta n]} < n,\; \BB_n\rcr\\
& \le & K \p\lcr A_{2^m} < 1 + 2^m + 4^{m-n\gamma} \;\;\forall \; m = \sigma_1 + 1,\ldots, \sigma_{[\delta n]} + 1, \;\; \sigma_{[\delta n]} < n\rcr.
\end{eqnarray*}
Introducing the notations
$$\GG_{[\delta n]}\; =\; \FF_{\sigma_{[\delta n]}} \quad\mbox{and}\quad \DD_{[\delta n]}\; =\; \lacc A_{2^m} < 1 + 2^m + 4^{m-n\gamma} \;\;\forall \; m = \sigma_1 + 1,\ldots, \sigma_{[\delta n]} + 1\racc$$
for every $n\ge 0,$ the strong Markov property at time $2^{\sigma_{[\delta n]}}$ yields
\begin{eqnarray*}
\p\lcr \DD_{[\delta n]}, \; \sigma_{[\delta n]} < n \rcr & \le & \p\lcr \Un^{}_{\DD_{[\delta n]-1}\cap\{ \sigma_{[\delta n]} < n\}} \p\lcr A_{2^{\sigma_{[\delta n]} + 1}} < 1 + 2^{\sigma_{[\delta n]} + 1} + 4^{\sigma_{[\delta n]} + 1-n\gamma}  \; \vert \; \GG_{[\delta n]}\rcr\rcr\\
& \le &  \p\lcr \Un^{}_{\DD_{[\delta n]-1}\cap\{ \sigma_{[\delta n]} < n\}} \p\lcr A_{2^{\sigma_{[\delta n]} + 1}} < 1 + 2^{\sigma_{[\delta n]} + 1} + 2^{2 +2\sigma_{[\delta n]}/\alpha}  \; \vert \; \GG_{[\delta n]}\rcr\rcr
\end{eqnarray*}
where in the second line we used the event $\{ \sigma_{[\delta n]} < n\}$. On the other hand, again by the strong Markov property, conditionally on $\GG_{[\delta n]}$ we can decompose
$$A_{2^{\sigma_{[\delta n]} + 1}}\; = \; A_{2^{\sigma_{[\delta n]}}} \; +\; 2^{\sigma_{[\delta n]}} Z_{2^{\sigma_{[\delta n]}}}\; +\;\hA_{2^{\sigma_{[\delta n]}}}$$ 
where $\hA$ is a copy of $A$ independent of $\GG_{[\delta n]}$. By the definition of $\sigma_{[\delta n]},$ this yields 
\begin{eqnarray*}
\p[ A_{2^{\sigma_{[\delta n]} + 1}} < 1 + 2^{\sigma_{[\delta n]} + 1} + 2^{2 +2\sigma_{[\delta n]}/\alpha} \, \vert \, \GG_{[\delta n]}] & \le & \p[ \hA_{2^{\sigma_{[\delta n]}}} < 2^{2 +\sigma_{[\delta n]}(1+ 1/\alpha)} + 2^{2 +2\sigma_{[\delta n]}/\alpha} \, \vert \, \GG_{[\delta n]}]\\
& \le & \p[ \hA_{2^{\sigma_{[\delta n]}}} < 2^{3 +\sigma_{[\delta n]}(1+ 1/\alpha)} \; \vert \; \GG_{[\delta n]}]\\
&\le & \p[ A_1 < 8],
\end{eqnarray*}
where we used the fact that $\alpha > 1$ in the second line, and the
$(1+1/\alpha)$-self-similarity of $\hA$ independent of $\GG_{[\delta
  n]}$ in the third. But from Proposition 3.4.1 in \cite{ST}, we know
that $A_1$ is a real $\alpha-$stable random variable whose L\'evy-Khintchine exponent $\Phi(\lambda) = - \log \e [\ee^{i\lambda A_1}]$ is given by
$$\Phi(\lambda) \; =\;(\alpha +1)^{-1}\vert\lambda\vert^\alpha(1-{\rm i}\beta{\rm sgn}(\lambda)\tan(\pi\alpha/2)), \quad \lambda\in\r,$$
where $\beta$ is the skewness parameter of $Z_1$. In particular, its positivity parameter is $\p[ A_1 > 0] = \rho > 0$ and its scaling parameter belongs to $[1/3,1/2)$ independent of $\alpha$. This clearly entails - see Section 1.6 in \cite{ST} - that there exists $\kappa < 1$ depending only on $\rho$ such that $\p [A_1 \ge 8] \ge (1-\kappa)$ for
every $\alpha > 1.$ Let me also stress that this argument breaks down when there are no negative jumps, since then $\p [A_1 \ge 8] \to 0$ when $\alpha \to 1.$ We finally obtain
\begin{eqnarray*}
\p\lcr \DD_{[\delta n]}, \; \sigma_{[\delta n]} < n \rcr & \le & \kappa \p\lcr \DD_{[\delta n]-1}, \; \sigma_{[\delta n]} < n \rcr\;\le \; \kappa \p\lcr \DD_{[\delta n]-1}, \; \sigma_{[\delta n]
  -1} < n -1 \rcr\;\le\;\kappa^{[\delta n]}
\end{eqnarray*}
by an induction argument, and since above the event $\{\sigma_{[\delta n]} <
n\}$ is only used to obtain the upper bound $4^{\sigma_{[\delta n]} +
  1-n\gamma}\le 2^{2 +2\sigma_{[\delta n]}/\alpha}$. Now because $\delta, \kappa$ are independent of $\alpha$, this yields (\ref{main2}) as desired, and completes the proof of Theorem B when $\alpha >1$. The case $\alpha \le 1$ can be handled exactly in the same way in working on the events
$$\lacc A_{2^m} < 1 + 2^m, \; m = 0 \ldots n\racc, \quad n \ge 0,$$
and we leave the details to the reader. However, for the sake of completeness and since our arguments are partly different from \cite{CD1}, we will give the \\

\noindent
{\bf Proof of the lemma.} Set $\p_n [.]= \p [\,. \,\vert \BB_n]$ for concision. By the
definition of $\sigma_{[\delta_n]},$ we have for every $\delta >0$
\begin{eqnarray*}
\p_n [ \sigma_{[\delta n]} < n] & = & 1 - \p_n\lcr \sum_{k=1}^n
\Un_{\CC_k^c} > \delta' n\rcr \; \ge \; 1 - \frac{1}{\delta' n} \sum_{k=1}^n \p[\CC_k^c]\;= \;  \frac{1}{\delta' n}
\sum_{k=1}^n \p_n[\CC_k]\; -\; \frac{\delta}{\delta'} 
\end{eqnarray*}
where we set $\delta' = 1-\delta$ for conciseness. Hence, it
suffices to prove the existence of $c > 0$ independent of
$\alpha$ and $k$ such that
$$\p_n[\CC_k] \; \ge \; c$$
for every $k\in [1, n]$ and every $n$ sufficiently large. To do so, we will first consider the events 
$$\DD_{k,n}\; =\; \lacc Z_t < 2^{-(n-k)/\alpha},\;\, \forall\, t\in [0,2^{n+k}]\racc\; \subseteq\; \BB_n$$
and prove the intuitively obvious inequalities
\begin{equation}
\label{holley}
\p[\CC_k\, \vert\, \BB_n] \; \ge \; \p[ \CC_k\, \vert\, \DD_{k,n}]
\end{equation}
for every $n\ge 1$ and every $k\in [1, n],$ with the help of a discretization of $Z$ and a FKG argument. Fixing $k$ and $n,$ set
$$Z^N_t\; =\; \sum_{i <tN} \lpa L_{\frac{i+1}{N}} - L_{\frac{i}{N}}\rpa\quad\mbox{and}\quad A^N_t\; =\;\frac{1}{N} \sum_{i <tN} \lpa \sum_{j\le i} \lpa L_{\frac{j+1}{N}} - L_{\frac{j}{N}}\rpa\rpa$$
for every $t \ge 0, N\ge 1.$ It is well-known - see e.g. Chapter 3 in \cite{Bil} - that the bivariate process $\{ (Z^N_t, A^N_t), \; t\in [0, 2^{n+k}]\}$ 
converges in law towards $\lacc (Z_t, A_t), \; t\in [0, 2^{n+k}]\racc$ for the Skorokhod topology when $N\to\infty$. Hence setting $\BB_n^N,\CC_k^N, \DD^N_{k,n}$ for the events $\BB_n,\CC_k, \DD_{k,n}$ with $(Z, A)$ replaced by $(Z^N, A^N)$, by weak convergence and right continuity it suffices to show
\begin{equation}
\label{holly}\p[\CC_k^N\, \vert\, \BB_n^N] \; \ge \; \p[ \CC_k^N\, \vert\, \DD_{k,n}^N]
\end{equation}
for every $N\ge 1$. The key-point is that the probabilities of the events $\BB_n^N,\CC_k^N, \DD^N_{k,n}$ and their intersections depend only on the joint law $\p^N$ of the increments $\{ (Z_{\frac{i}{N}} - Z_{\frac{i-1}{N}}), \; i = 1, \ldots, N2^{n+k}\},$ which are stationary and independent. The measure $\p^N$ has the density
\begin{equation}
\label{prod}
f(x) \; = \; \prod_{i=1}^{N2^{n+k}} g(x_i)
\end{equation}
with respect to the Lebesgue measure on $\r^{N2^{n+k}},$ where $g$ is the density of the variable $Z_{\frac{1}{N}},$
and the density of the conditional measure $\p_n^N[.] = \p^N [\, .\,\vert {\hat \BB}^N_n]$ with respect to the Lebesgue measure on $\r^{N2^{n+k}}$ is hence given by
$$h(x) \; = \;\frac{\Un_{{\hat \BB}_n^N} f(x)}{\p[\BB_n^N]}$$
where ${\hat \BB}_n^N$ is the image of $\BB_n^N$ under  $\p^N$. Besides, it is obvious from the definition of $\BB_n^N$ that the function $\Un_{{\hat \BB}_n^N}(x)$ is decreasing, in the sense that if $x_i \ge y_i$ for all $i = 1\ldots N2^{n+k},$ then $\Un_{{\hat \BB}_n^N}(x)\le \Un_{{\hat \BB}_n^N}(y).$ From (\ref{prod}) and this monotonicity property, we deduce that $h$ satisfies Holley's criterion:
$$h(x)h(y)\; \ge\; h(x\vee y)h(x\wedge y), \qquad x,y \in\r^{N2^{n+k}}$$
with the notations $(x\vee y)_i = x_i\vee y_i$ and $(x\wedge y)_i = x_i\wedge y_i$ for all $i = 1\ldots N2^{n+k}.$ By Corollary 12 in \cite{Ho}, this entails that the measure $\p_n^N$ satisfies FKG inequality in the sense that
$$\p_n^N[\CC\cap\DD] \; \ge \; \p_n^N[\CC]\p_n^N[\DD]$$
whenever $\Un_\CC$ and $\Un_\DD$ are increasing functions. On the other hand, the function $\Un_{{\hat \CC}_k^N}(x)$ is increasing and the 
function $\Un_{{\hat \DD}_n^N}(x)$ is decreasing, with the same notation as above for ${\hat \CC}_k^N$ and ${\hat \DD}_n^N.$ This finally entails
$$\p_n^N[{\hat \CC}_k^N\cap {\hat \DD}_n^N] \; \le \; \p_n^N[{\hat \CC}_k^N]\p_n^N[ {\hat \DD}_n^N]$$
and, after some elementary transformations using the inclusion $\DD_{k,n}^N\subseteq\BB_n^N$, the desired inequality (\ref{holly}) for every $k =1, \ldots, n$ and $n, N\ge 1$. Hence, from (\ref{holley}), we now need to show that there exists $c >0$ independent of $\alpha$ such that
$$\p[ \CC_k\, \vert\, \DD_{k,n}]\; \ge \; c$$
foe every $k\in[1,n]$ and $n$ sufficiently large. A scaling argument yields first
$$\p[ \CC_k\, \vert\, \DD_{k,n}]\; = \;{\hat \p}[Z_1 < 1, A_1 < 1\,\vert\, \{ Z_t > -2^{-n}, \; \forall t\in [0, 2^n]\}]$$
where ${\hat \p}$ stands for the law of the dual process $\hZ = -Z$. By Chaumont's results - see Remark 1 and Theorem 6 in \cite{Chaume}, the conditional law on the right-hand side converges to ${\hat \p}^{\uparrow}$ which is the law of $\hZ$ conditioned to stay positive. Hence, for $n$ sufficiently large, one has
$$\p[ \CC_k\, \vert\, \DD_{k,n}]\; \ge \;\frac{1}{2}\,{\hat \p}^{\uparrow}[Z_1 < 1, A_1 < 1]\; \ge\; \frac{1}{2}\,{\hat \p}^{\uparrow}[S_1 < 1]$$
with the notation $S_1 = \sup\{ Z_t, \, t \le 1\}$. It is now intuitively obvious that the right-hand side can be bounded from below by a positive constant on every closed interval $[\alpha_0, 2]\subset (0,2].$ Let us give a rigorous argument. Using Theorem 1 in \cite{Chomedu}, one gets
$${\hat \p}^{\uparrow}[S_1 < 1] \; \ge\; c\, {\hat \e}^{(me)}[Z^2_1\Un_{\{ S_1 < 1\}}]$$
for some constant $c$ depending only on $\rho$, where ${\hat \p}^{(me)}$ denotes the law of the meander associated to $Z$. The pathwise representation of the meander given in \cite{Ber}, Proposition VIII.16 yields 
$${\hat \e}^{(me)}[Z^2_1\Un_{\{ S_1 < 1\}}] \;\ge\; c'\,{\hat \e}[(Z_1 -I_1)^2 \Un_{\{S_1 -I_1< 1\}}]$$
for some constant $c'$ depending only on $\rho$, with the notation $I_1 = \inf\{ Z_t, \, t \le 1\}$. Finally, the L\'evy-Khintchine formula and Theorem 1 in \cite{Kada} entail that the probability on the right-hand side is a continuous function of the parameter $\alpha$. Since it is obviously positive for every $\alpha \in (0,2]$, it is bounded from below by a positive constant on every closed interval $[\alpha_0, 2]\subset (0,2]$, which completes the proof.\qed

\begin{REM} {\em The proof of the above lemma would be slightly simpler if we could prove that the measures $\p_n$ themselves satisfy FKG inequality, instead of considering their discretizations $\p^N_n$. From Example 2.3.6 and Theorem 4.6.1 in \cite{ST}, and by right-continuity of the sample paths of $Z$, we know that the unconditioned measure $\p$ satisfies FKG, in the sense for every time-horizont $T >0$ and every bounded measurable increasing functionals $F,G : \ddd_T \to \r^+,$
\begin{equation}
\label{FKG}
\e[F(Z_t, \, t\le T)G(Z_t, \, t\le T)]\; \ge\;\e[F(Z_t, \, t\le T)]\e[G(Z_t, \, t\le T)] 
\end{equation}
(here we set $\ddd_T$ for the Skorokhod space of c\`adl\`ag functions from $[0,T]$ to $\r$, and we say that a functional $F: \ddd_T\to\r^+$ is increasing if for every $x, y\in\ddd_T,$
$x_t\,\ge\, y_t \; \forall \, t\le T\; \Longrightarrow\; F(x_t, \, t\le T)\,\ge\, F(y_t, \, t\le T).$) Since $\BB_n$ is a monotonous event, our desired conditioned version of (\ref{FKG})
\begin{equation}
\label{FKG1}\e_n[F(Z_t, \, t\le T)G(Z_t, \, t\le T)]\; \ge\;\e_n[F(Z_t, \, t\le T)]\e_n[G(Z_t, \, t\le T)]
\end{equation} 
would be fulfilled if $\p$ satisfied the {\em strong FKG inequality.} However, there exist some path-measures which are FKG but not strong FKG - we learned this from J.-D. Deuschel, and we could not prove (\ref{FKG1}) directly, even in the Brownian case.}
\end{REM}

\section{Proof of Theorem A}
Recall that we are interested in the Hausdorff dimension of the random set
$$\LLA\; =\; \lacc a(x), \, x\in \r\;\;\mbox{and}\;\; a(x-)\, =\, a(x)\racc$$
where  $a(x) := \max\{ s \ge 0, \;C'(s) \le x\}$ and $C'$ is the right-derivative of $C$, the convex hull of the function 
$$x\;\mapsto\; \int_0^x (X_u + u) \, \d u,\qquad x\in\r.$$
Recall also that  $X$ is a two-sided $\alpha$-stable L\'evy process ($\alpha >1$) with positive jumps as defined in (\ref{XX}), and fix its positivity parameter $\rho < 1/\alpha$ once and for all. 

Notice that by definition $X$ does not jump negatively at time  $l$ whenever $l\in\LLA$: $X_l\ge X_{l-}$ a.s. On the other hand, it is well possible that $l\in\LLA$ is a "conical point" in the sense that $X_l > X_{l-}.$ However, if we define 
$${\hat \LLA}\; =\; \lacc a(x), \, x\in \r,\;\; a(x-)\, =\, a(x)\;\;\mbox{and}\;\; X_{a(x)} = X_{a(x)-}\racc,$$
we see from the fact that the set of points of discontinuity of $X$ is a.s. countable that
$${\rm Dim}_H \,\LLA\; =\; {\rm Dim}_H \,{\hat \LLA}\qquad \mbox {a.s.}$$
The key-point - which was first noticed by Sinai \cite{Si1} in the Brownian case - is that a.s.
$${\hat \LLA}\,\subseteq\, {\bar \LLA}\, :=\, \lacc a\in\r\; / \; \int_0^x  (X_u + u) \, \d u\, \ge\, \int_0^a  (X_u + u) \, \d u\, + \, (x-a) (X_a + a),\;\forall\; x\in \r\racc,$$
so that we only need to get an upper bound on ${\rm Dim}_H \,{\bar \LLA}$. To do so, we will use (\ref{main}) together with the same arguments as Molchan and Kholkhov \cite{MK}. First, a slight modification of Lemma 1 in \cite{MK} shows that Theorem A will be proved as soon as there exists a constant $\kappa > 0$ independent of $\alpha$ and a subsequence $\delta_n \to 0$ such that
$$\p\lcr {\bar \LLA} \cap (x-\delta_n, x+\delta_n)\neq \emptyset\rcr\; \le\; \delta_n^{\kappa}, \qquad n\to\infty$$
uniformly in $x\in\r$. Indeed, reasoning exactly as In Lemma 1 in \cite{MK} entails then 
$${\rm Dim}_H \,{\bar \LLA}\; \le\; 1 -\kappa\qquad \mbox{a.s.}$$
and completes the proof of the theorem with $\alpha_0 = 1/(1-\kappa).$ Second, we remark that by linearity of the integral and by independence and stationarity of the increments of the process $u\mapsto X_u +u,$ the random sets
$${\bar \LLA} \cap (x-\delta_n, x+\delta_n), \quad x\in\r$$
have all the same law. Hence, we need to show that there exists a constant $\kappa > 0$ independent of $\alpha$ such that
\begin{equation}
\label{Main1}
\liminf_{\delta\to 0}\delta^{-\kappa}\p\lcr {\bar \LL_{\delta}} \neq \emptyset\rcr\; = \; 0,
\end{equation}
where for every $\delta >0$ we wrote
$${\bar \LL_{\delta}}\; =\; \lacc \vert a\vert < \delta \;/ \; \int_0^x  (X_u + u) \, \d u\, \ge\, \int_0^a  (X_u + u) \, \d u\, + \, (x-a) (X_a + a),\;\forall\; x\in \r\racc.$$
For every $\delta >0,$ set $\FF_\delta$ for the completed filtration generated by $\{ X_u, \; \vert u\vert \le \delta\},$
and consider the $\FF_\delta$-measurable random variables
$$M_{\delta} \; =\; \sup_{\vert x\vert \le\delta}\vert X_x +x\vert\quad\mbox{and}\quad N_{\delta}\; =\;\sup_{\vert x\vert \le\delta}\lva \int_0^x (X_u +u)\, \d u\rva.$$
If $\hX = -X$ and $\hZ = -Z$ denote the dual processes of $X$ and $Z$ respectively, we see that a.s.
\begin{eqnarray*}
\lacc{\bar \LL_{\delta}} \neq \emptyset\racc &\subseteq &\lacc \int_0^x  (\hX_u - u) \, \d u\, \le\, N_\delta + (\delta + \vert x\vert)M_\delta\;\;\forall\; x\in \r\racc\;\;\subseteq\;\;\AA_\delta\,\cap\, \AA_\delta'
\end{eqnarray*}
with the notations
\begin{eqnarray*}
& \AA_\delta & =\;\; \lacc \int_\delta^x  ((\hZ_u - \hZ_\delta)- (u-\delta)) \, \d u\, \le\, 2 N_\delta + 2(\delta + x )M_\delta\;\;\forall\; x\ge \delta\racc \\
\mbox{and} & & \\
& \AA_\delta' & = \;\; \lacc \int_{-\delta}^{-x}  ((-\hZ_u' + \hZ_\delta')- (u-\delta)) \, \d u\, \le\, 2 N_\delta + 2 (\delta + \vert x\vert)M_\delta\;\;\forall\; x\ge \delta\racc,
\end{eqnarray*}
where $\hZ'$ is an independent copy of $\hZ$. By the Markov property, the events $\AA_\delta$ and $\AA_\delta'$ are independent and identically distributed conditionally on $\FF_\delta,$ so that
\begin{equation}
\label{AA}
\p\lcr {\bar \LL_{\delta}} \neq \emptyset\rcr\; =\; \p\lcr\p\lcr {\bar \LL_{\delta}} \neq \emptyset\;\vert\;\FF_\delta\rcr\rcr\; \le\; \e\lcr \p \lcr \AA_\delta \cap\AA_\delta'\; \vert\; \FF_\delta\rcr\rcr\; =\; \e\lcr \p \lcr \AA_\delta \; \vert\; \FF_\delta\rcr^2\rcr.
\end{equation}
Besides, again by the Markov property, we can write
$$\p \lcr \AA_\delta \; \vert\; \FF_\delta\rcr\; =\; \p\lcr \int_0^t L_u\, \d u \; \le \; 2n + 2 t m + t^2/2,\;\;\forall\, t \ge0\rcr_{n =N_\delta, m = M_\delta}$$
where $L$ is a copy of $\hZ$ independent of $\FF_\delta$. Now since $\alpha > 1$ and by the scaling property of $X$, it is immediate
to see that a.s. $N_\delta \le \delta^{1+1/\alpha} (1 + N)$ and $M_\delta \le \delta^{1/\alpha} (1 + M)$ as soon as $\delta <1,$ where $M, N$ are $\FF_\delta$-measurable and such that
$$M\; \elaw\;\sup_{\vert x\vert \le 1}\vert X_x\vert\quad\mbox{and}\quad N\; \elaw\;\sup_{\vert x\vert \le 1}\lva \int_0^x X_u\, \d u\rva.$$  
Setting $\varepsilon = (2\delta)^{(\alpha+1)/2\alpha}$ and $R = \max\{ (1+N), (1+M)^{\alpha +1}\}$ for conciseness, we can rewrite
\begin{eqnarray*}
\p \lcr \AA_\delta \; \vert\; \FF_\delta\rcr & = & \p\lcr \int_0^t L_u\, \d u \; \le \; \eps^2 r  +  (\eps^2 r)^{1/(\alpha +1)}t + t^2/2,\;\;\forall\, t \in [0,1]\rcr_{r =R}\cdot
\end{eqnarray*}
Returning to (\ref{AA}), we obtain
\begin{eqnarray*} 
\p\lcr {\bar \LL_{\delta}} \neq \emptyset\rcr & \le & \p\lcr R \ge \eps^{-1}\rcr\; +\; \lpa\p\lcr \int_0^t L_u\, \d u \; \le \; \eps  +  \eps^{1/(\alpha +1)}t + t^2/2,\;\;\forall\, t \in [0,1]\rcr\rpa^2\\
& \le & (c\delta)^{1/4} \; +\; \lpa\p\lcr \int_0^t L_u\, \d u \; \le \; \eps  +  \eps^{1/(\alpha +1)}t + t^2,\;\;\forall\, t \in [0,1]\rcr\rpa^2
\end{eqnarray*}
for some $c >0$, where in the third line we used well-known estimates on the upper tails of supremum of stable processes - see e.g. Theorem 10.5.1. in \cite{ST}. A scaling argument yields
$$\p\lcr \int_0^t L_u\, \d u \le  \eps  +  \eps^{1/(\alpha +1)}t + t^2,\;\forall\, t \in [0,1]\rcr\; =\; \p\lcr \int_0^t L_u\, \d u  \le  1  +  t + n^{-\gamma} t^2,\;\forall \,t \in [0,n]\rcr$$
where $n =\eps^{-\alpha/(\alpha +1)} = 1/\sqrt{2\delta}.$
Finally, since $L \elaw -Z$ has negative jumps, we see from (\ref{main}) that there exists $\kappa > 0$ depending only on $\rho$ such that 
$$\liminf_{\delta\to 0} \delta^{-\kappa}\p\lcr {\bar \LL_{\delta}} \neq \emptyset\rcr \; =\; 0,$$
which is (\ref{Main1}) and completes the proof.
\qed

\section{Two conjectures}

Let us start by a classical result of Bingham concerning the asymptotics of the ruin probabilities related to the stable process $Z$: if we set $S =\inf\{ t >0, \; Z_t >1\}$, then there exists a constant $c\in (0,\infty)$ such that
\begin{equation}
\label{bing}
\p[ S >t]\; \sim\; c t^{-\rho}, \qquad t\to +\infty
\end{equation}
as soon as $\vert Z\vert$ is not a subordinator - see Proposition VIII. 2 in \cite{Ber}. Notice in passing that the constant $c$ equals $\alpha \, p_{Z_1}(0)$ when $Z$ has no positive jumps - this is a consequence of Skorokhod's formula written e.g. in \cite{Bi} p. 749, and that in the other cases it can be given (non explicitly) in terms of the excursion measure associated to the reflected process - see Lemma 1 in \cite{Chomedu}. Considering now $T = \inf\{ t >0, \, A_t = 1\}$ the first-passage time of the integral of $Z$ across 1, from (\ref{bing}) it is tantalizing to state the 

\begin{CONC} Suppose that $\vert Z\vert$ is not a subordinator. Then
$$\p[ T >t]\; =\; t^{-\rho/2 +o(1)}, \qquad t\to +\infty.$$
\end{CONC}

The Brownian case $\alpha = 2$ and $\rho = 1/2$ had been obtained in \cite{JJG} after expanding a closed formula of McKean concerning the distribution of $T$, which yields actually a more accurate estimate like (\ref{bing}) with an explicit constant for $c$ - see Proposition 2 therein. The case with no negative jumps $\alpha > 1$ and $\rho = (\alpha -1)/\alpha$ was proved recently in \cite{ThS}, with a good control on $o(1)$ allowing to show that $\e[T^{\rho/2}] = +\infty.$ The above conjecture is also motivated by the aforementioned fact that $Z$ and $A$ have the same positivity parameter $\rho$, a quantity which should typically play a r\^ole in the distribution of $T$. On the other hand, the process $A$ is smoother than $Z$, so that the upper tails of the distribution of $T$ should be heavier than those of $S$. We propose the value $\rho/2$ for the critical exponent, since it is in accordance with the spectrally positive case. Conjecture C will be the matter of further research \cite{DSS}.

Suppose now $\alpha > 1$ and consider the drifted stable process $Z^c_t = Z_t + ct, \; t\ge 0,$ for some $c\neq 0$. From Lemma VI.21 in \cite{Ber} and explicit estimates on the renewal function of $Z^c,$ one can show that
$$\p[Z^c_t \, <\,\eps, \;\forall\; t\in [0,1]]\;\asymp\; \eps^{\rho\alpha},\qquad \eps\to 0$$
for every $c\in\r.$ Notice also that the latter estimate is false when $\alpha =1$ by Bingham's result, since $Z^c$ is then a strictly stable process whose positivity parameter depends on $c$. We think that the estimate is also untrue when $\alpha <1,$ but we got stuck in proving this. Setting $A^c$ for the integral of $Z^c$, we believe from the above fact that when $\alpha >1,$
\begin{equation}
\label{area}
\p[A^c_t \, <\,\eps, \;\forall\; t\in [0,1]]\; =\; \eps^{\rho\alpha/(\alpha +1) + o(1)},\qquad \eps\to 0
\end{equation}
for every $c\in\r.$ By self-similarity, notice that this estimate is the same as Conjecture C when $c =0$. When $c <0,$ it is particularly relevant for the inviscid Burgers equation whose initial data is the two-sided dual process $\hZ$. Indeed, reasoning as in the proof of Theorem A, one can show that the upper bound in (\ref{area}) entails ${\rm Dim}_H \, \LLA\; \le \; {\hat\rho}\;$ a.s. where ${\hat\rho} = 1-\rho$ is the positivity parameter of $\hZ$. Considering now (\ref{Burg}) where the initial data is the two-sided dual process $Z,$ from Conjecture C, the above considerations and Bertoin's result \cite{Ber1}, one is tempted to state the
\begin{COND} With the above notations, for every $\alpha > 1$ and every $\rho\in [1-1/\alpha, 1/\alpha]$ one has
 $${\rm Dim}_H \, \LLA\; = \; \rho\quad\mbox{a.s.}$$
\end{COND}
From the present paper, we are convinced that optimal lower tail estimates for the integral of $\hZ$ should provide the key-argument to obtain the upper bound in Conjecture D. The lower bound seems more delicate because of the positive jumps of $\hZ$ which prevent from using Handa's criterion \cite{panda}. On the other hand, it is immediate to see from its definition that the set $\LL_t$ contains the points of global increase of the drifted process $Z^{1/t}$ for every $t >0,$ and one may wonder if the methods developped by Marsalle \cite{marsala} to determine the Hausdorff dimension of the points of local increase of the non-drifted process $Z$, could not be useful. 

\bigskip 
\noindent {\bf Acknowledgement.} I am very grateful to F. Caravenna and J.-D.~Deuschel for having let me experience, live in Berlin, the ongoing of their paper \cite{CD1}.


\begin{thebibliography}{10}

\bibitem{Ber}
{\sc J.~Bertoin.} {\it L\'evy Processes.} Cambridge University Press, Cambridge, 1996.

\bibitem{Ber1}
{\sc J.~Bertoin.} The inviscid Burgers equation with Brownian initial velocity. {\em Commun. Math. Phys.} {\bf 193} (2), pp. 397-406, 1998.

\bibitem{Bil}
{\sc P.~Billingsley.} {\it Convergence of probability measures. 2nd ed.} Wiley, Chichester, 1999.

\bibitem{Bi}
{\sc N.~H.~Bingham.} Fluctuation theory in continuous time.
{\it Adv. Appl. Probability} {\bf 7} (4), pp. 705-766, 1975.

\bibitem{CD1}
{\sc F.~Caravenna} and {\sc J.-D.~Deuschel.} Pinning and wetting transition for (1+1)-dimensional fields with Laplacian 
interaction. Preprint, 2007.

\bibitem{CD2}
{\sc L.~Carraro} and {\sc J. Duchon.} Equation de Burgers avec conditions initiales \`a accroissements ind\'ependants et homog\`enes. {\em Ann. Inst. H. Poincar\'e Anal. Non Lin\'eaire} {\bf 15} (4), pp. 431-458, 1998.

\bibitem{Chaume}
{\sc L.~Chaumont.} Conditionings and path decomposition for L\'evy processes. {\em Stochastic Process. Appl.} {\bf 64} (1), pp. 39-54, 1996.

\bibitem{Chomedu}
{\sc L.~Chaumont.} Excursion normalis\'ee, m\'eandre et pont pour les processus de L\'evy stables. {\em Bull. Sci. Math.} {\bf 121} (5), pp. 377-403, 1997.

\bibitem{DSS}
{\sc A.~Devulder, Z. Shi} and {\sc T.~Simon.} The lower tail problem for integrated stable processes. Working paper.

\bibitem{JJG}
{\sc M.~Goldman.} On the first passage of the integrated Wiener process. {\em Ann. Math. Stat.} {\bf 42} (6), pp. 2150-2155, 1971.

\bibitem{panda}
{\sc K.~Handa.} A remark on shocks in inviscid Burgers' turbulence. In: Nonlinear waves and weak turbulence with applications in oceanography and condensed matter physics, pp. 339-345, {\em Progr. Nonlinear Differential Equations Appl.} {\bf 11}, Birkh\"auser, Boston, 1993.

\bibitem{Ho}
{\sc R.~Holley.} Remarks on the FKG inequalities. {\em Commun. Math. Phys.} {\bf 36} (2), pp. 227-231, 1974.

\bibitem{JW}
{\sc A.~W.~Janicki} and {\sc W.~A.~Woyczynski.}
Hausdorff dimension of regular points in stochastic Burgers flows with L\'evy $\alpha$-stable initial data. {\em J. Stat. Phys.} {\bf 86} (1-2), pp. 277-299, 1997.

\bibitem{Kada}
{\sc T.~V.~Kadankova.} On the joint distribution of the supremum, infimum, and the value of a semicontinuous process with independent increments. {\em Theor. Probab. Math. Statist.} {\bf 70}, pp. 61-70, 2005.

\bibitem{marsala}
{\sc L.~Marsalle.} Hausdorff measures and capacities for increase times of stable processes. {\em Potential Anal.} {\bf 9} (2), pp. 181-200, 1998. 

\bibitem{MK}
{\sc G.~Molchan} and {\sc A.~Khokhlov.} 
    Small values of the maximum for the integral of 
    fractional Brownian motion. {\it J. Stat. Phys.} {\bf 114} (3-4), pp. 923-946, 2004.
    
\bibitem{ST}
{\sc G.~Samorodnitsky} and {\sc M.~S.~Taqqu.} {\em Stable Non-Gaussian Random Processes.} Chapman \& Hall, New York, 1994.

\bibitem{ThS1}
{\sc T.~Simon.} The lower tail problem for the area of a symmetric
stable process. Talk given at the conference {\em Small Deviations Probabilities and Related Topics II}, 2005. Abstract available at {\tt http://www.pdmi.ras.ru/EIMI/2005/sd/talk/simon.pdf}

\bibitem{ThS}
{\sc T.~Simon.} The lower tail problem for homogeneous functionals of stable processes with no negative jumps. Preprint, 2007. Available at {\tt http://arxiv.org/pdf/math.PR/0701653}

\bibitem{Si1}
{\sc Ya.~G.~Sinai.} Statistics of shocks in solutions of inviscid Burgers equation. {\it Comm. Math. Phys.} {\bf 148} (3), pp. 601-621, 1992.
 
\bibitem{Wo}
{\sc W.~A.~Woyczynski.} {\em Burgers-KPZ turbulence. 
G\"ottingen lectures.} Lect. Notes Math. {\bf 1700} Springer-Verlag, Berlin, 1998.

\end{thebibliography}
\end{document}